\documentclass[12pt,leqno]{amsart}
\usepackage{amssymb}

\DeclareMathOperator*{\Ess}{\mathrm{ess\;sup}}
\DeclareMathOperator*{\ess}{\mathrm{ess\;inf}}

\newcommand{\al}{\alpha}
\newcommand{\barF}{\overline{F}}
\newcommand{\ds}{\displaystyle}
\newcommand{\e}{\varepsilon}
\newcommand{\eop}{\hfill$\square$}

\theoremstyle{plain}
\newtheorem{Thm}{Theorem}
\newtheorem{Cor}{Corollary}
\newtheorem{Lem}{Lemma}

\numberwithin{equation}{section}

\begin{document}

\title[The Mean Residual Life]{Limiting Behaviour of the Mean Residual
Life}

\date{Submitted March 7, 2001; revised May 29, 2002; accepted June 12, 2002}
\maketitle

\begin{center}
   {\sc David M. Bradley and Ramesh C. Gupta}

   {\em Department of Mathematics and Statistics, University of
        Maine, Orono, ME 04469-5752, U.S.A.}

\vskip.1in

   {\em e-mail: {\tt bradley}@{\tt math.umaine.edu},
                {\tt  rcgupta}@{\tt maine.maine.edu} }
\end{center}

\vskip.5in

%\author{David~M. Bradley}
%\address{Department of Mathematics \& Statistics\\
%         University of Maine\\
%         5752 Neville Hall
%         Orono, Maine 04469-5752\\
%         U.S.A.}
%\email[David~M. Bradley]{bradley@math.umaine.edu}
%\author{Ramesh~C. Gupta}
%\email[Ramesh~C. Gupta]{rcgupta@maine.maine.edu}

%\subjclass{Primary: 62N02; Secondary: 62N05, 62E20}
% 62N02 = Survival Analysis and censored data - estimation
% 62N05 =  "         "           "        "  - reliability and life testing
% 62E20 = Asymptotic Distribution Theory


%\keywords{Mean residual life, limiting behaviour, asymptotic
%expansion, failure rate, hazard function.}

\vskip.5in

%\begin{abstract}
{\sc Abstract.}
  In survival or reliability studies, the mean residual life or
  life expectancy is an important characteristic of the model.
  Here, we study the limiting behaviour of the
  mean residual life, and derive an asymptotic expansion which can
  be used to obtain a good approximation for large
  values of the time variable.  The asymptotic expansion
  is valid for a quite general class of failure rate
  distributions---perhaps the largest class that can be expected
  given that the terms depend only on the failure rate
  and its derivatives.
%\end{abstract}

\maketitle

\section{Introduction and Background}\label{sect:Intro}
In life testing situations, the expected additional lifetime given
that a component has survived until time $t$ is a function of $t$,
called the mean residual life.  More specifically, if the random
variable $X$ represents the life of a component, then the mean
residual life is given by $m(t)=E(X-t|X>t)$.  It is well known
that the mean residual life is related to the survival
(reliability) function $\barF$ by
\begin{equation}
   m(t) =\frac{1}{\barF(t)}\int_t^\infty \barF(u)\,du,
\label{MRLbarF}
\end{equation}
and to the failure rate (hazard function) $r=-{\barF}\,'/\,\barF$
by
\begin{equation}
   m(t)=\int_t^\infty \exp\bigg\{-\int_t^u r(x)\,dx\bigg\}\,du
   = e^{R(t)}\int_t^\infty e^{-R(u)}\,du,
\label{MRLR}
\end{equation}
where
\begin{equation}
   R(t) = \int_0^t r(x)\,dx = -\log \barF(t)
\label{Rdef}
\end{equation}
is the integrated failure rate (cumulative hazard function).  We
also have
\begin{equation}
   m'(t)=r(t)m(t)-1,
\label{DE1}
\end{equation}
(see Calabria and Pulcini (1987), for example.)

The mean residual life has been employed in life length studies by
various authors, e.g.~Hollander and Proschan (1975), Bryson and
Siddiqui (1969), and Muth (1977).  Muth (1977) observed that the
failure rate takes only the instantaneous present into account,
whereas the mean residual life takes the complete future into
account.  Meilijson (1972) has studied certain limiting properties
of the mean residual life.  A smooth estimator of the mean
residual life is given by Chaubey and Sen (1999).

In this paper, we undertake a detailed study of the limiting
behaviour of the mean residual life (\S\ref{sect:Limit}), and
derive an asymptotic expansion (\S\ref{sect:Asym}) which can be
used to obtain good approximations for large values of the time
variable. The asymptotic expansion is valid for a quite general
class of failure rate distributions---perhaps the largest class
that can be expected given that the terms depend only on the
failure rate and its derivatives.

For the family of age smooth distributions, Rojo (1996) has
established, for large values of the time variable, a relationship
between the mean residual life and the failure rate in terms of
the index $\rho$ of regular variation.  Calabria and Pulcini
(1987) noted a relationship between the limiting behaviour of the
mean residual life and the failure rate.   This relationship was
developed further by Chaubey and Sen (1999) for the class of
distributions having nondecreasing failure rate.  Our approach
provides a considerable improvement on their approximation, and
moreover does not require that the failure rate be nondecreasing.
%Of course, certain assumptions on the failure rate are required in
%order to have a valid asymptotic expansion; we would argue that
%our conditions are about as general and natural as one could hope
%for.

\section{Limiting Behaviour}\label{sect:Limit}
By applying L'H\^opital's rule to~(\ref{MRLbarF}), Calabria and
Pulcini~(1987) derived the relationship
\begin{equation}
   \lim_{t\to\infty}m(t)=\lim_{t\to\infty}\frac1{r(t)},
\label{CalPul}
\end{equation}
provided the latter limit exists and is finite.  They then
used~(\ref{DE1}) to conclude that $\lim_{t\to\infty} m'(t)=0$, or
equivalently, that
\begin{equation}
   \lim_{t\to\infty}r(t)m(t)=1.
   \label{reciprocal}
\end{equation}
Unfortunately, one cannot infer~(\ref{reciprocal})
from~(\ref{CalPul}) unless one also assumes that
$\lim_{t\to\infty} r(t)$ is finite and strictly positive.  For a
counterexample, fix positive constants $a$ and $b$ and consider
the linear mean residual life $m(t)=a+bt$ with corresponding
failure rate $r(t)=(1+b)/(a+bt)$.  The class of distributions with
linear mean residual life has been studied by Hall and
Wellner~(1981) and Oakes and Dasu~(1990).  Counterexamples that
satisfy $\lim_{t\to\infty} r(t)=\infty$ also exist, but as these
tend to be somewhat more complicated, further discussion is
deferred to the end of this section.  A more detailed study of the
limiting behaviour of the mean residual life follows.

\begin{Thm}\label{thm:ess}  Let $r$ denote the failure rate, and
let $m$ denote the corresponding mean residual life.   For $t>0$
define $\beta(t)=\Ess\{r(x):x>t\}$ and
$\alpha(t)=\ess\{r(x):x>t\}$. Then
\[
  \lim_{t\to\infty}\frac{1}{\beta(t)}\le
  \liminf_{t\to\infty}m(t)\le \limsup_{t\to\infty}m(t)
  \le \lim_{t\to\infty}\frac{1}{\alpha(t)}.
\]
\end{Thm}

\noindent{\bf Proof.}  Since $\alpha$ (respectively, $\beta$) is
clearly nondecreasing (nonincreasing),
\[
   m(t) = \int_0^\infty \exp\bigg\{ -\int_t^{x+t}
   r(y)\,dy\bigg\}\,dx
\]
implies
\[
   \int_0^\infty e^{-x\beta(t)}\,dx
   \le m(t)
   \le \int_0^\infty e^{-x\alpha(t)}\,dx.
\]
\eop

\begin{Cor}\label{Cor:Ess}  The failure rate $r$
and the mean residual life $m$ have limiting behaviour related by
\[
   \frac{1}{\ds\limsup_{t\to\infty}r(t)} \le \liminf_{t\to\infty}m(t)
   \le \limsup_{t\to\infty}m(t)\le
   \frac{1}{\ds\liminf_{t\to\infty}r(t)}.
\]
\end{Cor}

Note that Corollary~\ref{Cor:Ess} implies~(\ref{CalPul}) without
the assumption that $\lim_{t\to\infty}1/r(t)$ exists; for example,
if $\limsup_{t\to\infty}r(t)=0$, then
$\lim_{t\to\infty}m(t)=\infty$.  Corollary~\ref{Cor:Ess} also
implies that if $\lim_{t\to\infty}r(t)$ exists (finite) and is
strictly positive, then~\textup{(\ref{reciprocal})} holds.

The limiting reciprocal relationship~(\ref{reciprocal}) between
the failure rate and the mean residual life may be interpreted as
an approximation or asymptotic formula $m\sim s$, where $s=1/r$.
By imposing suitable conditions on $s$ and its derivatives, it is
possible to refine this approximation by introducing additional
terms into the asymptotic formula.  We shall carry out this
programme in \S\ref{sect:Asym}.  Even in cases where the
reciprocal relationship~(\ref{reciprocal}) fails, one can
sometimes obtain reasonably precise information about the limiting
behaviour of the product of the failure rate and the mean residual
life by studying the limiting behaviour of $s$ and its
derivatives.  By putting $E=[0,\infty)$ in
Corollary~\ref{Cor:Almost} below, it can be seen that for
continuously differentiable $s$, if $\lim_{t\to\infty}s'(t)=0$,
then~(\ref{reciprocal}) holds. However, it may happen that
$\lim_{t\to\infty} s'(t)$ does not exist.  It turns out that if
$|s'(t)|$ is not ultimately ``too big too often,'' one can still
say a good deal about the product $r(t)m(t)$ when $t$ is large:
see Theorem~\ref{Thm:EssSup} and Corollary~\ref{Cor:Almost} below.
There are also failure rates for which nonzero values of
$\lim_{t\to\infty}s'(t)$ are possible. In such cases, the
reciprocal relationship~(\ref{reciprocal}) may fail. Nevertheless,
we have the following result.
\begin{Thm}\label{Thm:LHop1} Denote the failure rate and its reciprocal by
$r$ and $s=1/r$, respectively, and let $m$ denote the
corresponding mean residual life.  Let $R$ denote the integrated
failure rate~\textup{(\ref{Rdef})}. Suppose that
\begin{equation}
   \lim_{t\to\infty} s(t)\exp(-R(t)) =0, \label{hyp1}
\end{equation}
and that
\begin{equation}
   \lim_{t\to\infty} (1-s'(t))^{-1}
   \mbox{\; \textup{exists (finite)}}.\label{hyp2}
\end{equation}
Then
\[
   \lim_{t\to\infty}r(t)m(t)
      = \lim_{t\to\infty} \left(1-s'(t)\right)^{-1}.
\]
\end{Thm}

\noindent{\bf Proof.}  Apply L'H\^opital's rule to
\[
   r(t)m(t)
      =\int_t^\infty e^{-R(x)}\,dx \bigg/\big(s(t)e^{-R(t)}\big).
\]
\eop

When applying Theorem~\ref{Thm:LHop1} in practice, it may often be
easier to verify conditions that imply the
hypothesis~(\ref{hyp1}), as opposed to verifying~(\ref{hyp1})
directly. We give an example of such a condition below.
\begin{Thm}\label{Thm:Hyp1} Denote the reciprocal of the failure rate
by $s$, and let $R$ denote the integrated failure
rate~\textup{(\ref{Rdef})}. Suppose that $\lim_{t\to\infty} \Ess\{
s'(x): x>t\} <1$ and that $s(t)$ is finite for all sufficiently
large values of $t$.  Then~\textup{(\ref{hyp1})} holds.
\end{Thm}

\noindent{\bf Proof.}  There exist numbers $L<1$ and $t_0\ge0$
such that for all $t\ge t_0$, $s(t)$ is finite and
$\Ess\{s'(x):x>t\}<L$.  For $t\ge t_0$, we have
\begin{align*}
   s(t)\exp(-R(t)) &= \exp\bigg\{\log s(t_0)+
         \int_{t_0}^t \frac{s'(x)}{s(x)}\,dx
         -\int_0^t r(x)\,dx \bigg\}\\
   &= s(t_0)\exp\bigg\{\int_{t_0}^t\big(s'(x)-1\big)r(x)\,dx
       -\int_0^{t_0}r(x)\,dx\bigg\}\\
   &\le  s(t_0)\exp\bigg\{\big(\Ess_{x>t_0} s'(x)-1\big)\int_{t_0}^t
   r(x)\,dx -R(t_0)\bigg\}\\
   &\le s(t_0)\exp\big\{(L-1)(R(t)-R(t_0))-R(t_0)\big\}\\
   &= s(t_0)e^{-LR(t_0)}\big(\barF(t)\big)^{1-L}.\\
\end{align*}
Since $\lim_{t\to\infty}\barF(t)=0$ and $L<1$,~(\ref{hyp1})
follows. \eop

A condition that implies the hypothesis $\lim_{t\to\infty} \Ess\{
s'(x): x>t\} <1$ of Theorem~\ref{Thm:Hyp1} and is typically even
easier to verify is $\limsup_{t\to\infty} s'(t)<1$.  Thus,
Theorems~\ref{Thm:LHop1} and~\ref{Thm:Hyp1} can be readily applied
to a wide variety of situations in which the long term behaviour
of the rate of change of the reciprocal of the failure rate is
known. A very simple yet illustrative example is provided by the
distribution with fractional linear failure rate $r(t)=1/(c+dt)$.
Here $s'$ has constant value $d$, and thus if $d<1$, then
Theorems~\ref{Thm:LHop1} and~\ref{Thm:Hyp1} imply that
$\lim_{t\to\infty}r(t)m(t)=1/(1-d)$.  Of course, in this case it
is easy to verify this fact directly, as the mean residual life
function is linear: $m(t)=(c+dt)/(1-d)$.  However,
Theorems~\ref{Thm:LHop1} and~\ref{Thm:Hyp1} are equally applicable
to cases in which it may be difficult or impossible to determine
the long term behaviour of the mean residual life directly.

The issue of practical usefulness apart, there are compelling
theoretical reasons which point up the significance of
Theorem~\ref{Thm:LHop1} as well.  Recall the observation
that~(\ref{reciprocal}) may be interpreted as an approximation or
asymptotic formula $m\sim s$. Viewing the conclusion of
Theorem~\ref{Thm:LHop1} in the same light yields
\[
   m\sim \frac{s}{1-s'} = s(1+s'+(s')^2+(s')^3+\cdots), \qquad |s'|<1,
\]
which agrees with the first three terms of the asymptotic
expansion~(\ref{asympt}).  Thus, in some sense,
Theorem~\ref{Thm:LHop1} represents a strengthening of Calabria and
Pulcini's first order result~(\ref{reciprocal})---subject to the
appropriate conditions---to third order.

We next address the problem of what happens when $s'$ exists but
$\lim_{t\to\infty}s'(t)$ does not.  As might be expected, in
general some information about the limiting behaviour of the
product $r(t)m(t)$ is lost. Nevertheless, in many cases one can at
least bound $rm$ in terms of the essential supremum of $s'$. We'll
see as a result that even if $\limsup_{t\to\infty} |s'(t)|>0$, as
long as $|s'(t)|$ is not ultimately ``too big too often''
then~(\ref{reciprocal}) holds.

\begin{Thm}\label{Thm:EssSup}  Let $r$ denote the failure rate, and
let $m$ denote the corresponding mean residual life.  Let $s=1/r$,
and put $\lambda(t) := \Ess\{|s'(x)| : x>t\}$.  Suppose that
$r(t)$ is positive and continuously differentiable for all
sufficiently large values of $t$, and that
$\lim_{t\to\infty}\lambda(t)<1$.  Then for all sufficiently large
values of $t$, $1/(1+\lambda(t))\le r(t)m(t)\le 1/(1-\lambda(t))$.
\end{Thm}

\noindent{\bf Proof.}  By Theorem~\ref{Thm:Hyp1},
condition~(\ref{hyp1}) holds.  Thus, the hypotheses permit us to
integrate by parts and discard the limit at the upper range of
integration in the integrated term.  We have for all sufficiently
large values of $t$,
\[
   r(t)m(t) = r(t)e^{R(t)}\int_t^{\infty}r(x)e^{-R(x)}s(x)\,dx
   = 1+r(t)e^{R(t)}\int_t^{\infty} e^{-R(x)}s'(x)\,dx.
\]
But
\[
   \bigg|\int_t^{\infty}e^{-R(x)}s'(x)\,dx\bigg|
   \le \Ess_{x>t}|s'(x)| \int_t^{\infty} e^{-R(x)}\,dx.
\]
Thus, for all sufficiently large values of $t$,
\[
   1-\lambda(t) r(t)m(t)\le r(t)m(t)\le 1+\lambda(t) r(t)m(t).
\]
\eop

It follows that the reciprocal relationship~(\ref{reciprocal}) may
hold even if $\limsup_{t\to\infty}|s'(t)|>0$, as long as
$\lim_{n\to\infty} s'(t_n)=0$ for ``most'' sequences
$t_1<t_2<\dots\to\infty$.  Here, ``most'' is in the sense of
Lebesgue measure.  The following result makes this observation
more precise.

\begin{Cor}\label{Cor:Almost}
Let $s$ denote the reciprocal of the failure rate.  Suppose that
$s(t)$ is finite and continuously differentiable for all
sufficiently large values of $t$.  Suppose further that there
exists a subset $E$ of the interval $[0,\infty)$ whose complement
in $[0,\infty)$ is of Lebesgue measure zero, and such that for
every sequence $t_1,t_2,\ldots$ of elements of $E$ with
$\lim_{n\to\infty} t_n=\infty$ we have
$\lim_{n\to\infty}s'(t_n)=0$. Then~\textup{(\ref{reciprocal})}
holds.
\end{Cor}

\noindent{\bf Proof.}  As customary, denote the indicator function
of a set $A$ by $\chi_A$.  Suppose that $E\subseteq [0,\infty)$
satisfies the hypotheses of the Corollary. Then the complement of
$E$ in $[0,\infty)$ has Lebesgue measure zero, and
$\lim_{x\to\infty}s'(x)\chi_{E}(x)=0$. Therefore, if $\e>0$ is
given, there exists a suitably large value of $t$ such that
$|s'(x)\chi_{[t,\infty)}(x)|<\e$ for all $x\in E$. By definition
of the essential supremum, this implies (in the notation of
Theorem~\ref{Thm:EssSup}) that $\lambda(t)=\Ess\{|s'(x)| : x>t \}
<\e$. As $\lambda$ is a nonincreasing function and $\e>0$ is
arbitrary, it follows that $\lim_{t\to\infty}\lambda(t)=0$. \eop

We conclude this section with an example of a distribution in
which $\lim_{t\to\infty} r(t)=\infty$, but
$\lim_{t\to\infty}r(t)m(t)\ne 1$.  Let $a,b,c,d$ be positive
constants satisfying $a>b$, $d>2b$, and $c^2>(a+b)d$. Consider the
mean residual life defined by
\[
   m(t) = \frac{a+b\sin(t^2)}{c+dt},\qquad t\ge 0.
\]
We have
\[
   m'(t) =
   \frac{2bt\cos(t^2)}{c+dt}-\frac{(a+b\sin(t^2))d}{(c+dt)^2}.
\]
Thus,
\[
   \liminf_{t\to\infty} m'(t) = -2b/d,\qquad \limsup_{t\to\infty}
   m'(t)=2b/d,
\]
and hence
\[
   \liminf_{t\to\infty}r(t)m(t)=1+\liminf_{t\to\infty}m'(t)=1-2b/d<1,
\]
whereas
\[
   \limsup_{t\to\infty}r(t)m(t)=1+\limsup_{t\to\infty}m'(t)=1+2b/d>1.
\]
Since
\begin{align*}
   (c+dt)^2(1+m'(t)) &=
   (c+dt)^2+(c+dt)2bt\cos(t^2)-(a+b\sin(t^2))d\\
   &\ge c^2+2cdt+d^2t^2-2bct-2bdt^2-(a+b)d\\
   &= (d-2b)dt^2+(d-b)2ct+c^2-(a+b)d
\end{align*}
is clearly positive for $t\ge 0$, it follows that
$r(t)=(1+m'(t))/m(t)>0$ for $t\ge 0$.  Furthermore,
\[
   r(t)=\frac{c+(d+2b\cos(t^2))t}{a+b\sin(t^2)}-\frac{d}{c+dt}
\]
implies $\lim_{t\to\infty}r(t)=\infty$.

\section{Asymptotic Expansion}\label{sect:Asym}

Under certain conditions, the mean residual life has an asymptotic
expansion in terms of the failure rate and its derivatives.  An
initial attempt in this direction was made by Chaubey and Sen
(1999) for the class of distributions having nondecreasing failure
rate.  However, it is easy to envision situations where, say with
regular maintenance, even an ultimately nondecreasing failure rate
may be an inappropriate model.  Therefore, we provide an
alternative approach that requires no monotonicity assumptions on
the failure rate.

We take Chaubey and Sen's asymptotic formula
\[
   m(t)=\frac{1}{r(t)}-\frac{r'(t)}{(r(t))^3}
   +O\bigg(\frac{r''(t)}{(r(t))^4}\bigg),\qquad t\to\infty,
\]
(with error term corrected), as a point of departure.
Unfortunately, their derivation is not rigorous, and as a result,
certain growth conditions on $r$ and its derivatives are omitted.
Nevertheless, their approach can, in principle, be made to work,
and thus one can show that under suitable conditions on the
failure rate there is an asymptotic series development to
arbitrary order that begins
\[
   m\sim
   r^{-1}-r'r^{-3}-r''r^{-4}+\left(3(r')^2-r'''\right)r^{-5}
   +\left(10r'r''-r''''\right)r^{-6}+\cdots,\qquad t\to\infty.
\]
Here, we have made the abbreviations $m=m(t)$, $r^{-1}=1/r(t)$,
etc.   More explicitly, for each positive integer $n$, we have
\begin{equation}
   m \sim \sum_{k=0}^{n-1} c_kr^{-k-1} + (\mbox{additional\,
   terms}),
\label{mAsym}
\end{equation}
where as $t\to\infty$, the additional terms tend to zero more
rapidly than $c_kr^{-k-1}$ for $1\le k\le n-1$. Here, $c_0=1$,
$c_1=0$, $c_2=-r'$, $c_3=-r''$, $c_4=3(r')^2-r'''$, and in
general, $c_k$ is a polynomial in $r',r'',\dots,r^{(k-1)}$ for
$k\ge 2$, the explicit form of which is given by
\begin{equation}
   c_k = c_k(t) = k! \sum_{p=0}^{\lfloor k/2\rfloor} (-1)^p \sum
   \prod_{j\ge 1}
   \frac{1}{\al_j!}\bigg(\frac{r^{(j)}(t)}{(j+1)!}\bigg)^{\al_j},
\label{ckDef}
\end{equation}
where the inner sum is over all nonnegative integers
$\al_1,\al_2,\ldots$ such that $\sum_{j\ge 1}\al_j=p$ and
$\sum_{j\ge 1}(j+1)\al_j=k$.

Rather than attempt a rigorous proof of~(\ref{mAsym})
and~(\ref{ckDef}) along these lines, we instead develop an
alternative, operator-theoretic approach that appears to be
simpler, yet more powerful.  By way of motivation, observe that
since we are interested in the situation when~(\ref{reciprocal})
holds, it makes sense to let $s=1/r$ and rewrite the differential
equation~(\ref{DE1}) in the form $(1-sD)m = s$, where $D$ denotes
the derivative operator. Abbreviating $sD$ by $\Theta$, we might
hope to find a meaningful way to define the inverse operator
$(1-\Theta)^{-1}$ in such a way that given suitable growth
conditions on $s$ and its derivatives, we have a legitimate
asymptotic expansion
\[
   m = (1-\Theta)^{-1} s \sim
   s+\Theta(s)+\Theta^2(s)+\Theta^3(s)+\cdots.
\]
That this is indeed the case is the main content of
Theorem~\ref{Thm:Asym1}, in which an explicit formula is also
given for $\Theta^k(s)$, the $k$th term of the expansion.  A
significant step in this direction is provided by the following
result.

\begin{Lem}\label{Lem:IP} Let $s$ denote the reciprocal of the failure
rate, and suppose that for some nonnegative integer $n$, $s$ is
$n+1$ times continuously differentiable on the positive real line.
As usual, let $R$ denote the integrated failure
rate~\textup{(\ref{Rdef})}, and let $D$ denote the derivative
operator. For positive integers $k$, let $\Theta^k: = (sD)^k$ and
for convenience, let $\Theta^0:=s$. If in addition, we have
$\lim_{x\to\infty} \exp(-R(x)) \Theta^k(s(x))=0$ for each
$k=0,1,\dots,n$, then the mean residual life can be expressed in
the form
\begin{equation}
   m(t) = s(t) + \sum_{k=1}^{n} \Theta^k(s(t)) + e^{R(t)}
   \int_t^\infty r(x)e^{-R(x)}\Theta^{n+1}(s(x))\,dx,\qquad t\ge 0.
\label{mIP}
\end{equation}
\end{Lem}

\noindent{\bf Proof.}  The conditions on $s$ permit us to
integrate by parts $n+1$ times and discard the limits at the upper
range of integration.  Thus, integrating by parts once, we have
\[
   m(t)= e^{R(t)}\int_t^\infty r(x)e^{-R(x)}s(x)\,dx
    = s(t) + e^{R(t)} \int_t^\infty r(x)e^{-R(x)}\Theta(s(x))\,dx.
\]
It is now evident that the claimed formula~(\ref{mIP}) can be
proved using mathematical induction and repeatedly integrating by
parts.

Alternatively, note that $(1-s(t)D) \exp(R(t))\int_t^\infty
r(x)\exp(-R(x))h(x)\,dx = h(t)$ for all locally integrable
functions $h$ for which the integral converges, and so we may {\em
define}
\[
   (1-\Theta)^{-1} h(t) = e^{R(t)}\int_t^\infty r(x)e^{-R(x)}
   h(x)\,dx
\]
for all such functions $h$.  If we now take $h=s$ and then
$h=\Theta^{n+1} s$, then in light of the identity
\[
   (1-\Theta)^{-1}=
   1+\Theta+\Theta^2+\cdots+\Theta^{n}+(1-\Theta)^{-1}\Theta^{n+1},
\]
we have
\begin{align*}
   m(t) &= e^{R(t)} \int_t^\infty r(x)e^{-R(x)}s(x)\,dx\\
   &= (1-\Theta)^{-1}s(t)\\
   &= s(t)+\sum_{k=1}^{n}\Theta^k(s(t))+
   (1-\Theta)^{-1}\Theta^{n+1}(s(t))\\
   &= s(t)+\sum_{k=1}^{n}\Theta^k(s(t))+e^{R(t)}\int_t^\infty
   r(x)e^{-R(x)}\Theta^{n+1}(s(x))\,dx,
\end{align*}
as claimed. \eop

\begin{Thm}\label{Thm:Asym1}  Let $s$ denote the reciprocal of the
failure rate.  Suppose that for some nonnegative integer $n$, $s$
is $n+1$ times continuously differentiable on the positive real
line, and that $\lim_{t\to\infty}s(t)\exp(-R(t))=0$, where as
usual $R$ denotes the integrated failure
rate~\textup{(\ref{Rdef})}. Let $D$ and $\Theta$ be as in
Lemma~\ref{Lem:IP}.  Suppose that for each $k=0,1,\dots,n$ we have
$s(t) s^{(k+1)}(t)=o(s^{(k)}(t))$ as $t\to\infty$. Then
$\Theta^{k+1}(s)=o(\Theta^k(s))$ for all $0\le k\le n$, and we
have the asymptotic expansion
\[
   m(t)=s(t)+\sum_{k=1}^{n}\Theta^k(s(t))
   + o\big(\Theta^{n}(s(t))\big),
   \qquad t\to\infty
\]
for the mean residual life.  The terms $\Theta^k(s)$ may be
expressed in terms of $s$ and its derivatives by means of
\begin{equation}
   \Theta^k(s)=s\sum_{j_1,\dots,j_k\ge 0} d(j_1,\dots,j_k)
   \prod_{p=1}^k s^{(j_p)},
\label{Theta}
\end{equation}
where the sum extends over all nonnegative integers $j_p$ $(1\le
p\le k)$ and the coefficients $d(j_1,\dots,j_k)$ are generated by
the polynomial equation
\begin{equation}
   \sum_{j_1,\dots,j_k\ge 0} d(j_1,\dots,j_k)\prod_{p=1}^k
   x_p^{j_p} = \prod_{p=1}^k \sum_{j=1}^p x_j
\label{dGF}
\end{equation}
in the indeterminates $x_1,\dots,x_k$.
\end{Thm}

\noindent{\bf Proof.} The stated formula for $\Theta^k(s)$ is a
special case of a more general result.  See Snyder (1982).  We
next show that, under the hypotheses on $s$,
$\Theta^{k+1}(s)=o(\Theta^k(s)))$ for $0\le k\le n$. Since
$\Theta^{k+1}(s)=sD\Theta^k(s)$, we have
\begin{align*}
   \Theta^{k+1}(s) &= sDs\sum_{j_1,\dots,j_k\ge0}d(j_1,\dots,j_k)
   \prod_{p=1}^k \left(D^{j_p}s\right)\\
   &= s'\Theta^k(s)+s\sum_{j_1,\dots,j_k\ge 0}d(j_1,\dots,j_k)
   sD\prod_{p=1}^k\left(D^{j_p}s\right)\\
   &= s'\Theta^k(s)+s\sum_{j_1,\dots,j_k\ge0}d(j_1,\dots,j_k)
   \sum_{q=1}^k\left(s D^{j_q+1}s\right)\prod_{\substack{p=1\\p\ne q}}^k
   \left(D^{j_p} s\right).
\end{align*}
Note that the generating function~(\ref{dGF}) implies that
$j_1+\cdots+j_k=k$ above.  In particular, since $0\le k\le n$,
each $j_q\le n$, and so each $sD^{j_q+1}s = o(D^{j_q}s)$.  Also,
$s'=o(1)$ because $ss'=sDs=o(s)$.  Therefore, we have
\begin{align*}
   \Theta^{k+1}(s) &= o\left(\Theta^k(s)\right)+o\bigg(\sum_{q=1}^k
   s\sum_{j_1,\dots,j_k\ge 0}d(j_1,\dots,j_k)\prod_{p=1}^k
   D^{j_p}s\bigg)\\
   &= o\left(\Theta^k(s)\right)+o\left(k\,\Theta^k(s)\right)\\
   &= o\left(\Theta^k(s)\right).
\end{align*}
Next, observe that since $\lim_{t\to\infty}s(t)\exp(-R(t))=0$ and
$\Theta^{k+1}(s)=o(\Theta^k(s))$ for all $k=0,1,\dots,n$, we have
$\lim_{t\to\infty}\exp(-R(t))\Theta^k(s(t))=0$ for all
$k=0,1,\dots,n$. Therefore, in light of Lemma~\ref{Lem:IP}, it
remains only to show that
\[
   m(t)-\sum_{k=0}^n \Theta^k(s(t))=
   e^{R(t)}\int_t^\infty r(x)e^{-R(x)}\Theta^{n+1}(s(x))\,dx
   = o\left(\Theta^n(s(t))\right),
   \qquad t\to\infty.
\]
Since $\Theta^{n+1}(s)=o(\Theta^n(s))$, Lemma~\ref{Lem:IP} gives
\begin{align*}
   m(t)-\sum_{k=0}^n\Theta^k(s(t))
   &= e^{R(t)}\int_t^\infty r(x)e^{-R(x)}\Theta^{n+1}(s(x))\,dx\\
   &= o\bigg(e^{R(t)}\int_t^\infty
   r(x)e^{-R(x)}\Theta^n(s(x))\,dx\bigg)\\
   &= o\bigg(m(t)-\sum_{k=0}^{n-1}\Theta^k(s(t))\bigg).
\end{align*}
Therefore,
\[
   (1-o(1))\bigg(m(t)-\sum_{k=0}^{n-1}\Theta^k(s(t))\bigg)
   = \Theta^n(s(t)),
\]
from which it follows that
\[
   m(t)-\sum_{k=0}^{n-1}\Theta^k(s(t))
   = \Theta^n(s(t))(1+o(1)),
\]
or equivalently,
\[
   m(t)-\sum_{k=0}^{n}\Theta^k(s(t))=
   o\left(\Theta^{n}(s(t))\right),
\]
as claimed. \eop

We have
$\Theta^0(s)=s$, $\Theta(s)=ss'$, $\Theta^2(s)=s(s')^2+s^2s''$,
$\Theta^3(s)=s(s')^3+4s^2s's''+s^3s'''$, and
$\Theta^4(s)=s(s')^4+11s^2(s')^2s''+4s^3(s'')^2+7s^3s's'''+s^4s''''.$
Thus, under the conditions of Theorem~\ref{Thm:Asym1}, our
asymptotic expansion begins
\begin{equation}
   m\sim s+ss'+s(s')^2+s^2 s''+s(s')^3+4s^2s's''+s^3s'''+\cdots,
   \qquad t\to\infty.
\label{asympt}
\end{equation}

\noindent{\bf Acknowledgement.}  The authors are grateful to the
anonymous referees for a careful reading of the original
manuscript.  Their suggestions contributed to improvements in the
exposition.

\end{document}